\documentclass[11pt, reqno]{amsart}
\usepackage{amsmath,amssymb}
\usepackage[latin1]{inputenc}
\usepackage{epsfig}
\usepackage[all]{xy}

\newtheorem{teo}{Theorem}[]
\newtheorem{cor}[teo]{Corollary}
\newtheorem{pro}[teo]{Proposition}
\newtheorem{lem}[teo]{Lemma}

\theoremstyle{remark}

\newtheorem{exe}[]{Example}

\begin{document}
\newcommand{\ii}{\'{\i}}
\renewcommand{\r}{\hspace*{1 cm}}
\newcommand{\codim}{\mbox{codim\,}}

\newcommand{\cre}{{\rm Cr}(\plan)}
\newcommand{\aut}{{\rm Aut}}
\newcommand{\jon}{{\rm Jon}}
\newcommand{\bir}{{\rm Bir}}
\newcommand{\rat}{{\rm Rat}}
\newcommand{\hir}{{\mathbb F}}
\newcommand{\pic}{{\rm Pic}}
\newcommand{\rank}{\rm rg\,}
\newcommand{\pgl}{{\rm PGL}}
\newcommand{\gl}{{\rm GL}(\integer)}
\newcommand{\gln}{{\rm GL}_n(\integer)}
\newcommand{\mn}{{\rm M}_n(\integer)}
\newcommand{\voln}{{\rm Vol}_n}
\newcommand{\vol}{{\rm Vol}}

\newcommand{\bl}{{\mathbb B}l}
\newcommand{\mor}{\mbox{Mor}}
\newcommand{\ara}{\ar@}

\renewcommand{\div}{{\rm div}}
\newcommand{\supp}{{\rm Supp}}
\newcommand{\ord}{{\rm ord}}
\newcommand{\hcd}{{\rm hcd}}
\renewcommand{\r}{\hspace*{1 cm}}
\newcommand{\ie}{c'est-\`a-dire}
\newcommand{\Ie}{C'est-\`a-dire}
\newcommand{\proj}{{\mathbb P}}
\newcommand{\projn}{{\mathbb P}^n}
\newcommand{\plan}{{\mathbb P}^2}
\newcommand{\complex}{{\mathbb C}}
\newcommand{\corp}{{\mathbb K}}
\newcommand{\torus}{{\mathbb T}}
\renewcommand{\natural}{{\mathbb N}}
\newcommand{\real}{{\mathbb R}}
\newcommand{\rational}{{\mathbb Q}}
\newcommand{\integer}{{\mathbb Z}}
\newcommand{\grass}{{\mathbb G}}
\newcommand{\poids}{\mbox{poids}}
\newcommand{\wtilde}{\widetilde}

\newcommand{\te}{\text}
\newcommand{\lto}{\longrightarrow}
\newcommand{\tor}{\xymatrix{\ar@{-->}[r]&}}

\newcommand{\cali}{{\mathcal I}}
\newcommand{\calx}{{\mathcal X}}
\newcommand{\cala}{{\mathcal A}}
\newcommand{\cale}{{\mathcal E}}
\newcommand{\calt}{{\mathcal T}}
\newcommand{\caln}{{\mathcal N}}
\newcommand{\calf}{{\mathcal F}}
\newcommand{\calg}{{\mathcal G}}
\newcommand{\calh}{{\mathcal H}}
\newcommand{\cald}{{\mathcal D}}
\newcommand{\calo}{{\mathcal O}}
\newcommand{\calw}{{\mathcal W}}
\newcommand{\calr}{{\mathcal R}}
\newcommand{\calp}{{\mathcal P}}
\newcommand{\calu}{{\mathcal U}}
\newcommand{\call}{{\mathcal L}}
\newcommand{\calm}{{\mathcal M}}
\hyphenation{In-ter-sec-tion  in-ter-sec-tion as-so-ci ral-ment
trans-for-ma-tion ge-ne-ri-cal-ly}
\title[Determinantal Cremona Transformations]{On Characteristic Classes of
  Determinantal Cremona Transformations}
\author{Gerard Gonzalez-Sprinberg \and Ivan Pan}
\address{Gerard Gonzalez-Sprinberg\\ Institut Fourier\\Université de Grenoble
  I\\38402 St. Martin d'Hères\\France}\email{gonsprin@ujf-grenoble.fr}
\address{Ivan
Pan\\ Instituto de Matem\'atica\\UFRGS\\av. Bento Gon\c calves
9500\\91540-000 Porto Alegre,RS\\Brasil}\email{pan@mat.ufrgs.br}


\begin{abstract}
We compute the multidegrees and the Segre numbers of general
determinantal Cremona transformations, with generically reduced
base scheme, by specializing to the standard Cremona
transformation and computing its Segre class via mixed volumes of
rational polytopes.
\end{abstract}

\maketitle

\textsl{Dedicated to the memory of Shiing-Shen Chern}

\section{Introduction}

Let $F$ be a  birational map of degree $n$ of the projective space
$\proj^n$, given by the maximal minors of a $(n+1)\times n$ matrix
with general linear forms as entries, over an algebraically closed
field $\corp$. Such a \emph{determinantal Cremona transformation}
may be defined geometrically by $n$ correlations in general
position and they have been considered in the classical literature
on Cremona transformations by several authors, e. g. \cite{Cay},
\cite{Ber} and \cite[Chap. VIII, \S 4]{SR}.

The family of base schemes of determinantal maps, not necessarily
birational, may be identified as an open and connected subset of
the Hilbert scheme of the 2-codimensional arithmetically
Cohen-Macaulay subschemes of $\proj^n$ (see \cite{Ell}).

In this article we compute the \emph{multidegrees} of such $F$ and
the \emph{Segre classes} of its base scheme $B_F$, by specializing
to the \emph{standard Cremona transformation}
\[S_n:=(X_1X_2\cdots X_n:\cdots:X_0\cdots\widehat{X_i}\cdots
X_n:\cdots:X_0X_1\cdots X_{n-1})\] and by applying methods of
toric geometry. In this way it translates into computing mixed
volumes of some special polytopes with integer vertices.

The sequence of multidegrees $(d_0,\ldots,d_k,\ldots,d_n)$,
classically called the ``type'' of a Cremona transformation
$T:\proj^n\tor\proj^n$, are given by the degrees of the (direct)
strict transforms by $T$ of general $k$-dimensional linear
subvarieties $H^{n-k}$ of $\proj^n$; for a reference on this
``type'' via intersection theory see \cite{ST}. The multidegrees
are closely related (Prop. \ref{pro5}) to the \emph{Segre class}
of $T$,  $s(B_T,\proj^n)$, defined as the inverse of the Chern
class of the normal bundle of the embedding $i:
B_T\hookrightarrow\proj^n$ of the base scheme $B_T$ of $T$, if
$B_T$ is regularly embedded. The general definitions are given
later. This Segre class lives in the Chow group
$A_{\centerdot}(B_T)$ of $B_T$ ; we also consider its image
$s(B_T)$ in $A_\centerdot(\proj^n)$ and the \emph{Segre numbers}
$s_k = \int s(B_T) \cdot H^k$.

The main results are the following:

\begin{teo}
The determinantal Cremona transformations of $\, \proj^n$ with
generically reduced base scheme and the standard Cremona
transformation $S_n$ have the same multidegrees and Segre numbers.
\label{teo1}
\end{teo}

\begin{teo}
Let $n$ be an integer, $n\geq 2$. The multidegrees of the standard
Cremona transformation $S_n$ are equal to the binomial
coefficients:
\[d_k={n\choose k},\; 0\leq k\leq n. \]
\label{teo2}
\end{teo}

\begin{teo}
The image of the Segre class of $S_n$ in the Chow group
$A_{\centerdot}(\proj^n)$ is
\[s(B):=i_* s(B,\proj^n)=\sum_{k=0}^{n-2} s_k[H^{n-k}]\]
with Segre numbers
\begin{eqnarray*}
s_{k}&=&(-1)^{n-k-1}\int_{\proj^n} \pi_*(E^{n-k})\cdot
H^{k}\\
&=&(-1)^{n-k-1}\sum_{j=0}^{n-k} (-1)^{j}{n-k\choose j}{n\choose j} n^{n-k-j}
\end{eqnarray*}
for $0\leq k\leq n-2$, where $\pi:\bl_{B}(\proj^n)\to \proj^n$ is
the blowing-up of $\proj^n$ with center the base scheme
$B=B_{S_{n}}$ and where $E$ is the exceptional divisor of $\pi$.
\label{teo3}
\end{teo}

\begin{cor}
The Chow group of the base scheme $B$ of $S_n$ is
\[A_k(B)=\left\{\begin{array}{ll}\integer[H^{n-k}]_B&\mbox{for}\ 0\leq k\leq
    n-3\\
\bigoplus_{i=1}^{n(n+1)/2} \integer[\alpha_i]_B&\mbox{for}\
k=n-2\end{array}\right.
\]
where the $\alpha_i$ are the irreducible components of the reduced
base scheme $|B|$ of $S_n$, the $2$-codimensional skeleton of the
arrangement of the fundamental hyperplanes of $\proj^n$. With the
above notations, the Segre class of $S_n$ is
\[s(B,\proj^n)=\sum_{k=0}^{n-3}
s_k[H^{n-k}]_{B}+\left([\alpha_1]_B,\ldots,[\alpha_{n(n+1)/2}]_B\right)\]
where the $s_k$ are the Segre numbers of Theorem~\ref{teo3}.
\label{cor4}
\end{cor}

The multidegrees of the standard Cremona transformation $S_n$ are
well known for small dimensions as the ``types'' of $S_n$. We
didn't find a classical proof for higher dimensions. Some other
basic references related to (standard) Cremona transformations or
toric geometry, besides those appearing in the following text,
are: \cite{Dem}, \cite{DO}, \cite{KKMS}.

 We thank I.
Dolgachev for his advice on this subject, and F. Russo for useful
conversations on conservation of numbers for flat deformations. We
are pleased to acknowledge the support of Univ. Fed. RGS at Porto
Alegre, Institut Fourier of Univ. Grenoble 1 and the France-Brazil
Cooperation Program.

\section{Segre classes and multidegrees}

Let $F=(F_0:\cdots F_n):\proj^n\tor\proj^n$ be a rational map.
Denote by $B=B_F$ its base scheme, the subscheme of $\proj^n$ of
codimension at least two defined by the homogeneous ideal
$(F_0,\ldots,F_n)$ . Let $\pi:\bl_{B}(\proj^n)\to\proj^n$ be the
blowing-up of $\proj^n$ with center $B$, and $E:=\pi^{-1}(B)$ its
exceptional divisor. This blowing-up resolves the indeterminacies
of $F$ and we get a commutative diagram
\[\xymatrix{E\,\ar@{^(->}[r]\ar@{->}[d]_{\pi|_{E}}&\bl_B(\proj^n)\ar@{->}[d]_\pi\ar@{->}[rd]^\varphi& \\
B\,\ar@{^(->}[r]^{i}&\proj^n\ar@{-->}[r]^F& \projn}\]

The Segre class of $B$ in $\projn$ is defined as the Segre class
in the Chow group $A_\centerdot(B)$ of the normal cone
$C_B(\projn)$ ; it is equal to
\[s(B,\proj^n):=\sum_{i\geq 1}(-1)^{i-1} (\pi|_{E})_*(E^i)\]
(see~\cite[Cor. 4.2.2]{Ful}). We consider
its image in $A_\centerdot(\projn)$,
\[s(B):=i_*s(B,\projn)=\sum s_k[H^{n-k}]\]
where the integers $s_k=\int_{\projn}s(B)\cdot H^k$ are  called
the \emph{Segre numbers} of $F$. The Segre classes and numbers are
invariant under birational morphisms (\cite[Prop. 4.2]{Ful}).

The \emph{multidegrees} $d_k$ of $F$, $0\leq k\leq n$, defined as
the degrees of the (direct) strict transforms by $F$ of general
linear subvarieties $H^{n-k}$ of dimension $k$ of $\projn$, may be
computed via a resolution of indeterminacies of $F$.

 Let $H_1$ (resp. $H_2$) be the (total) pullback of a hyperplane $H$ by $\pi$
(resp. by $\varphi$). Then
\[d_k=d_k(F)=\int_{\bl_B(\projn)} H_1^{n-k}\cdot H_2^{k}.\]

Let $\partial=\deg(F):=\deg(F_i), \forall\,i$; then
$\varphi^*\calo(1)=\pi^*\calo(\partial)\otimes\calo(-E)$. It follows
$H_2\sim \partial H_1-E$.

\begin{pro}
The Segre numbers $s_k$ and the multidegrees $d_k$ of a rational map $F$ are related by:

(i) $d_0=1$, $d_1=\partial$, and $d_k=\partial^k-\sum_{0\leq \ell\leq k} {k\choose
  \ell}\partial^{k-\ell} s_{n-\ell}$ for $2\leq k\leq n$.

(ii) $s_k=(-1)^{n-k-1}\sum_{0\leq \ell\leq n-k} (-1)^\ell{n-h\choose \ell}\partial^{n-k-\ell}d_\ell.$
\label{pro5}
\end{pro}
\begin{proof}
(i) The equality follows from
\[d_k=\int_{\bl_B(\projn)} H_1^{n-k}\cdot H_2^k=\int_{\bl_B(\projn)}
H_1^{n-k}\cdot (\partial H_1-E)^k,\]
\[s_{n-\ell}=(-1)^{\ell-1}\int_{\projn} H_1^{n-\ell}\cdot\pi_*(E^\ell)\ ,\ \ \int_{\bl_B(\projn)} H_1^n=1\]
and the projection formula.\\
\ \\
(ii) From (i) we have
\[(d_0,\ldots,d_n)^t=M(-1,0,s_{n-2},s_{n-3},\ldots,s_0)^t,\]
where $M=(a_{k\ell})$ is the lower triangular matrix with entries
$a_{k\ell}=-{k\choose \ell}\partial^{k-\ell}$ if $0\leq \ell\leq
k\leq n$. It follows that the inverse matrix $M^{-1}=(b_{k\ell})$
is given by
$b_{k\ell}=(-1)^{1+k+\ell}{k\choose\ell}\partial^{k-\ell}$ if
$0\leq\ell\leq k\leq n$ and $0$ otherwise. It is enough to check
that $\chi(n,k):=\sum_{k\leq\ell\leq
  n} (-1)^\ell{n\choose\ell}{\ell\choose k}=0$ for $k<n$; this follows from
$\chi(n,k)=-\chi(n-1,k-1)$ by Stieffel formula and by induction on $k$ since
$\chi(n,0)=0$ for $n>0$.
\end{proof}

\section{General determinantal and standard Cremona transformations}

If $F$ is a determinantal Cremona transformation, according to
\cite{Ell} $B_F$ may be considered as a point of an open connected
and smooth set of the Hilbert scheme parameterizing the
arithmetically Cohen-Macaulay schemes of codimension 2, whose
projecting cones admit a determinantal resolution of length 2 (see
\cite[Thm. 2(ii)]{Ell}). Remark that for the base scheme $B=B_F$
under consideration there exists a minimal resolution
\begin{equation}
\xymatrix{0\ar@{->}[r]&\calo_{\proj^n}^{n}(-n-1)\ar@{->}[r]^{M_F}&\calo_{\proj^n}^{n+1}(-n)
\ar@{->}[r]&\cali_{B}\ar@{->}[r]&0}\label{res1}
\end{equation}
where $\cali_B$ denotes the ideal sheaf of $B$ and $M_F$ is the
$(n+1)\times n$ matrix of linear forms attached to $F$.

The standard Cremona transformation is indeed determinantal: it
may be defined by the maximal minors of the matrix
$M_{S_n}=(m_{ij})$ with entries given by
\[
m_{ij}=\left\{\begin{array}{l} \delta_{ij}X_{j-1}\ \ \mbox{if}\ \ i<n+1\\
-X_n\ \ \mbox{if} \ \ i=n+1, j=1,\ldots,n\end{array}\right.
\]

By the Peskine-Szpiro deformation Theorem~\cite[Thm. 6.2]{PS}
there exists a dense open set $S$ of a projective space over
$\corp$ and a $S$-scheme $X$ of codimension 2 in
$\proj^n_S=\proj^n\times S$, flat over $S$, such that

(a) the ideal sheaf $\cali_X$ of $X$ has a minimal resolution
\begin{equation}\xymatrix{0\ar@{->}[r]&\calo_{\proj^n_S}^n(-n-1)\ar@{->}[r]&\calo_{\proj^n_S}^{n+1}(-n)
\ar@{->}[r]&\cali_{X}\ar@{->}[r]&0}\label{res2}
\end{equation}

(b) for each codimensional 2 scheme $B$ having a minimal
resolution as in (\ref{res1}) there is a unique point $s\in S$
with $X(s)=B\times \{s\}$ such that the minimal resolution of the
associated ideal sheaf is obtained from (\ref{res2}) by
tensorizing with $\corp(s)$ over $\calo_S$.

\ \\
 \textsl{Proof of } \textbf{Theorem~\ref{teo1}} (stated in the
introduction):

\ \\
Suppose that $F$ is a general determinantal Cremona transformation
in the sense of the statement  ; with the above notations let $f$
and $s_n$ in $S$ be the points corresponding to the base schemes
of $F$ and $S_n$ respectively in the Hilbert scheme. Let $T$ be
the intersection of $S$ with the line joining $f$ and $s_n$.

For $s\in S$ let $B_s\subset\proj^n_s=\proj^ n\times\{s\}$ the
corresponding arithmetically Cohen-Macaulay scheme of codimension 2.

The set of $s\in T$ such that $B_s$ is
generically reduced contains a dense open set, since $B_f$ and $B_{s_n}$ are
generically reduced.

Let $\calx_T\subset\projn_T:=\proj^n\times T$ be the flat family
of the $B_t$'s parameterized by $t\in T$. Consider the blowing-up
$\pi:\widetilde{\projn_T}\to\projn_T$ of $\projn_T$ with center
$\calx_T$, let $E$ be the exceptional divisor; denote by
$\pi_t:\widetilde{\projn_t}\to \projn_t$ the corresponding
blowing-up in level $t\in T$. Choose general sections $\calh_1\in
\pi^*\calo_{\projn_T}(n)$ and $\calh_2\in
\pi^*\calo_{\projn_T}(n)\otimes \calo_{\widetilde{\projn_T}}(-E)$.

The family $(\calh_2)_t$ is flat over $T$ (see \cite[App. B.6.7]{Ful}). For
 general $t\in T$, since the divisor
 $\pi_t^*((\calh_1)_t)$ is regular in codimension 1 then it is normal. It
 follows that $(\calh_1)_t$ is a
 generically  flat family over $T$.

By the ``conservation of number'' (see \cite[Cor. 10.2.1]{Ful})
one has
\[(\calh_1)^k_f\cdot (\calh_2)^{n-k}_f=(\calh_1)^k_{s_n}\cdot
(\calh_2)^{n-k}_{s_n},\] so we get the equality of multidegrees,
and the equality of Segre numbers follows from the preceding
Proposition. \qed

\section{Fans and mixed volumes}

Consider the projective space $\projn$ over $\corp$ as the toric
variety defined by the fan $\Delta$ associated to the faces of the
simplex with vertices $e_0,e_1,\ldots,e_n$ where $e_1,\ldots,e_n$
is the standard basis of $\real^n$ and $e_0:=-\sum_{i=1}^n e_i$.

The central symmetry $-Id$ in $\real^n$ induces the standard
Cremona transformation $S_n$ viewed as a monomial birational map
(\cite{GSP}). A natural way to resolve the indeterminacies of
$S_n$ is to factorize it through the toric variety $X_{\Sigma}$
associated to the fan given by the minimal common subdivision of
$\Delta$ and $-\Delta$.
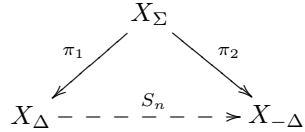
\[\xymatrix{&X_\Sigma\ar@{->}[ld]_{\pi_1}\ar@{->}[rd]^{\pi_2}&\\
X_\Delta\ar@{-->}[rr]^{S_n}& & X_{-\Delta}}\]

Intersection numbers for divisors in toric varieties may be computed as mixed
volumes of Minkowsky sums of polytopes (\cite{Oda}, \cite{Ful}), and we
obtain:

\ \\
\emph{Proof of} \textbf{Theorem~\ref{teo2}} (stated in the
introduction):

\ \\
 The multidegrees are independent of the resolution
of indeterminacies, so we may use the toric variety $X_\Sigma$.
Let $H_1$ (resp. $H_2$) be the pullback by $\pi_1$ (resp. by
$\pi_2$) of the hyperplane associated to the one-dimensional cone
through $e_0$ (resp. through $-e_0$). It turns out that the
polytope associated to $H_1$ (resp. to $H_2$) is the simplex
$\delta_n=[0,e_1,\ldots,e_n]$ (resp. the simplex $-\delta_n$). It
follows that the intersection number $d_k=\int_{X_\Sigma}
H_1^k\cdot H_2^{n-k}$ equals the mixed volume
$V(\delta_n,k;-\delta_n,n-k)$; i.e. the coefficient of the
monomial $\nu_1\cdots \nu_n$ in the homogeneous polynomial
$\vol\left((\nu_1+\cdots+\nu_k)\delta_n+(\nu_{k+1}+\cdots+\nu_n)(-\delta_n)\right)$.

\begin{lem}
Let $a,b$ be non negative integers. Then
\[\vol(a\delta_n+b(-\delta_n))=\sum_{j=0}^n {n\choose j} \frac{a^j b^{n-j}}{j
 !(n-j)!} \]
\label{lem6}
\end{lem}
\begin{proof}
This is a variant of the Steiner decomposition formula. It may be proved by
intersecting $a\delta_n+b(-\delta_n)$ with the $2^n$ orthants of
$\real^n$. Notice that each intersection is a cartesian product of simplices, with
${n\choose j}$ products of type $a\delta_j\times b(\delta_{n-j})$, for each $j,\
0\leq j\leq n$.
\end{proof}

Then the value of $d_k$ follows by putting $a=\nu_1+\cdots+\nu_k$,
$b=\nu_{k+1}+\cdots+\nu_n$. The only term with the monomial
$\nu_1\cdots\nu_n$ appears in $a^k b^{n-k}$ with coefficient $k!(n-k)!$, so
the lemma gives $d_k={n\choose k}$.

\qed

\section{Last proofs and examples}

The proof of \textbf{Theorem~\ref{teo3}} (stated in the
introduction)  follows from  Proposition~\ref{pro5},
Theorem~\ref{teo2} and the fact that $\deg(S_n)=n$.

\ \\

The reduced base scheme $|B|$ of $S_n$ is the union of the $n(n+1)/2$
codimension 2 linear subspaces obtained by cutting pairwise the fundamental
hyperplanes. This is a connected set for $n\geq 3$. By induction using the
fundamental exact sequence
\[A_\centerdot(X\cap Y)\to A_\centerdot(X)\oplus A_\centerdot(Y)\to
A_\centerdot(X\cup Y)\to 0\]
one gets that $A_k(B)$ is an infinite cyclic group, for $0\leq k\leq n-3$, and
$A_{n-2}(B)\simeq \integer^{n(n+1)/2}$ with a natural basis given by the
fundamental classes in $B$ of the irreducible components.

If $0\leq k\leq n-3$ then $i:B\hookrightarrow\projn$ induces an
isomorphism from $A_k(B)$ onto $A_k(\projn)$. By
Theorem~\ref{teo3}, the Segre number $s_{n-2}=n(n+1)/2,$ which is
the number of irreducible components of $|B|$. Then the symmetries
of $B$ give the $(n-2)$-dimensional part of $s(B,\projn)$, which
proves \textbf{Corollary~\ref{cor4}}(stated in the introduction).
\qed

\begin{exe}

A Cremona transformation $F:\proj^n\tor\proj^n$ is called
\emph{special} if its base scheme $B_F$ is smooth and irreducible.
It follows from \cite{ESB} (see also \cite{CK}) that if a
determinantal Cremona transformation is special then $n=3,4,5$ and
in this case $\deg(F)=\deg(F^{-1})=n$; moreover, if $n\neq 4$
these are the only special Cremona transformations. From
\cite[thm. 6.2]{PS} such a transformation exists. The equality on
degrees of $F$ and $F^{-1}$ is a particular case of our theorems
\ref{teo1} and \ref{teo2}.

In \cite{K} the case $n=3$ is considered (see also \cite[chap.
VIII, \S 4]{SR}). In particular the image of the Segre class of
$B_F$ in $\proj^3$ is computed there:
\[s(B_F)=B_F-28p\]
where $B_F$ is a smooth curve of degree 6 and genus 3 and $p$ is
the class of a point. This result confirms the number obtained for
$S_3$.
\end{exe}

\begin{exe}
Some examples of Segre numbers of a general determinantal Cremona
transformation or $S_n$, obtained following  theorems \ref{teo1}
and \ref{teo3} :

\[ s_{n-2}= \frac{1}{2} n(n+1) , \; s_{n-3}=-\frac{1}{3}
n(n+1)(2n+1), \; s_{n-4}=\frac{1}{8}
n(n+1)(5n^2+5n+2),\]
\[s_{n-5}=-\frac{1}{30} n(n+1)(2n+1)(2n^2+7n+6)\]
for $n\geq 7$; on the other hand $s_0$ for $n=2, 3, 4, 5$ is $3,
-28, 255, -2376$ respectively.

In general the formula given in theorem \ref{teo3} is equivalent
to an expression in terms of generalized hypergeometric function

\[s_k=-hypergeom([-n,k-n],[1],-\frac{1}{n})(-n)^{n-k}\].

The precise meaning of this expression may be found and verified
with \textsl{Maple}, by computing the formula for the Segre
numbers $s_k$ we obtained.

\end{exe}

\begin{exe}
Let $F:\proj^3\tor\proj^3$ be the Cremona transformation
\[F=(X_0(X_1X_2-X_0X_3):X_1(X_1X_2-X_0X_3):X_0X_1(X_2-X_3):X_0X_1(X_0-X_1));\]
its inverse is
\[F^{-1}=(X_0X_3:X_1X_3:X_0(X_2-X_1):X_1(X_2-X_0)).\]
This transformation $F$ is defined by the maximal minors of
\[M_F=\left(\begin{array}{rrc}
0&X_1&0\\
X_0&0&0\\
-X_1&-X_1&X_1-X_0\\
X_3&X_3&X_2-X_3\end{array}\right).\]
Remark that $d_2(F)=d_1(F^{-1})=2$, but for the standard Cremona
transformation of $\proj^3$ we have $d_1=d_2=3$. A direct verification shows
that $B_F$ is not reduced along its irreducible component $X_0=X_1=0$.

For higher $n$ we obtain analogous examples by induction: complete
the matrix $M_F$ with a line of zeroes and then a column  of general linear
forms to obtain a $4\times 5$ matrix. The maximal minors define a determinantal Cremona
transformation $G:\proj^4\tor\proj^4$. The matrix $M_{G^{-1}}$ can be obtained
from $M_{F^{-1}}$ in the same way that one obtains $M_{G}$ from $M_{F}$. We
obtain $\deg(G)>\deg(G^{-1})$.
\end{exe}

\begin{exe}
In $\proj^3$ there exists an involutory  determinantal Cremona
transformation whose base scheme is a sextic of arithmetic genus 3
supported on a (double) twisted cubic $\gamma$. Following
\cite[pag. 180]{SR}, such transformation may be constructed by
associating to a generic point $p\in\proj^3$ the intersection of
the polar planes of $p$ with respect to three general quadrics
containing $\gamma$. Theorem \ref{teo1} does not apply for this
type of determinantal transformation, but its Segre numbers are
also those of $S_3$.
\end{exe}


\begin{thebibliography}{99}
\bibitem{Ber} L. Berzolari, \emph{Algebraische Transformationen und
    Korrespondenzen}, Encyclopädie der Mathematischen Wissenschaften, dritter
    Band: Geometrie, 2.2.B. Teubner, 1932.
\bibitem{Cay} A. Cayley, \emph{On the rational Transformations in the two
    Spaces}, Proc. London Math. Soc., 3:127-180, 1870.
\bibitem{CK} B. Crauder, S. Katz, \emph{Cremona Transformations with smooth
  irreducible fundamental locus}, Amer. J. of Math., 111, 289-309, 1989.
\bibitem{Dem} M. Demazure, \emph{Sous-groupes alg\'ebriques de rang maximum du groupe de
    Cremona},
    Ann. Scient. Éc. Norm. Sup., $4^e$ série, 4,
    507-588, 1970.
\bibitem{DO} I. Dolgachev, D. Ortland, \emph{Point Sets in Projective Spaces and
    Theta Functions}, Ast\'erisque, 165, 1988.
\bibitem{Ell} G. Ellingsrud, \emph{Sur le schéma de Hilbert des variétés de
    codimension 2 dans $\proj^e$ à cône de
    Cohen-Macaulay}, Ann. Scient. Éc. Norm. Sup., $4^e$ série, 8,
    423-432, 1975.
\bibitem{ESB} L. Ein, N. Shepherd-Barron, \emph{Some Special Cremona
    Transformations}, Amer. J. Math. 111, 783-800, 1989.
\bibitem{Ful} W. Fulton, \emph{Intersection Theory}, Springer-Verlag, 1984.
\bibitem{Ful1} W. Fulton, \emph{Introduction to Toric Varieties}, Princeton
  University Press, 1993.
\bibitem{GSP} G. Gonzalez-Sprinberg, I. Pan, \emph{On the Monomial
    Birational Maps of the projective Space}, An. Acad. Bras. Ciências, 75(2), 129-134, 2003.
\bibitem{K} S. Katz, \emph{The Cubo-Cubic Transformation of $\proj^3$ is Very
    Special}, Math. Z. 195, 255-257, 1987.
\bibitem{KKMS} G. Kempf, F. Knudsen, D. Mumford, B. Saint-Donat \emph{Toroidal embeddings
    I}, LNM 339, Springer-Verlag, 1973.
\bibitem{Oda} T. Oda, \emph{Convex Bodies and Algebraic Geometry},
  Springer-Verlag, 1988.
\bibitem{PS} C. Peskine, L. Szpiro, \emph{Liaison des variétés
    algébriques. I}, Inventiones Math. Vol. 26, 271-302,
    1976.
\bibitem{SR} J. G. Semple, L. Roth, \emph{Introduction to Algebraic Geometry},
  Oxford at the Clarendon Press, 1949.
\bibitem{ST} J. G. Semple, J. A. Tyrrel, \emph{Specialization of Cremona
    Transformations}, Mathematika, 15, 171-177, 1968.
\end{thebibliography}
\end{document}